\documentclass{tran-l}
 \newtheorem{thm}{Theorem}[section]

 \newtheorem{prop}[thm]{Proposition}
 \theoremstyle{definition}
 
 \theoremstyle{remark}
 \newtheorem{rem}[thm]{Remark}
 \numberwithin{equation}{section}

\begin{document}

\title[Approximation of the Multiplication Table Function]
 {Approximation of the Multiplication Table Function}

\author{Mehdi Hassani}

\address{Institute for Advanced Studies in Basic Sciences\\
P.O. Box 45195-1159\\
Zanjan, Iran.}

\email{mmhassany@srttu.edu}

\thanks{}



\date{}

\dedicatory{}

\commby{}
\begin{abstract}
In this paper, considering the concept of Universal Multiplication
Table, we show that for every $n\geq 2$, the inequality:
$$
M(n)=\#\{ij|1\leq i,j\leq n\}\geq\frac{n^2}{\mathfrak{N}(n^2)},
$$
holds true with:
$$
\mathfrak{N}(n)=n^{\frac{\log 2}{\log\log
n}\left(1+\frac{387}{200\log\log n}\right)}.
$$
\end{abstract}

\maketitle

\textbf{Note.} In the first version of this paper, there are some
great mistakes, which I have done them refereing to a reference in
internet (comparing two versions, you can find that mistakes).
Professor Kevin Ford mentioned me that mistakes and announced my
some very interested improvements concerning the results of this
paper (see Remark \ref{ford-comments} at the end of this papre). I
deem my duty to thank him for his very kind comments.

\section{Introduction}
Consider the following $n\times n$ \textit{Multiplication Table},
which we denote it by $MT_{n\times n}$:
\begin{center}
\begin{tabular}{|c|c|c|c|c|}
  \hline
  1 & 2 & 3 & $\cdots$ & $n$ \\
  \hline
  2 & 4 & 6 & $\cdots$ & $2n$ \\
 \hline
  3 & 6 & 9 & $\cdots$ & $3n$ \\
 \hline
  $\vdots$ & $\vdots$ & $\vdots$ & $\ddots$ & $\vdots$ \\
 \hline
  $n$ & $2n$ & $3n$ & $\cdots$ & $n^2$ \\
  \hline
\end{tabular}
\end{center}
Let $\mathfrak{M}(n;k)$ be the number of $k$'s, which appear in
$MT_{n\times n}$; i.e.
\begin{equation}\label{mnk}
\mathfrak{M}(n;k)=\#\left\{(a,b)\in\mathbb{N}_n^2~|~ab=k\right\},
\end{equation}
where $\mathbb{N}_n=\mathbb{N}\cap [1,n]$. For example, we have:
$$
\mathfrak{M}(2;2)=2, \mathfrak{M}(7;6)=4, \mathfrak{M}(10;9)=3,
\mathfrak{M}(100;810)=10, \mathfrak{M}(100;9900)=2.
$$
In this paper first we study some elementary properties of the
function $\mathfrak{M}(n;k)$, for a fixed $n\in\mathbb{N}$. Then we
try to connect $\mathfrak{M}(n;k)$ by the famous
\textit{Multiplication Table Function}\footnote{This sequence, has
been indexed in ``The On-Line Encyclopedia of Integer Sequences''
data
base with ID A027424. Web page of above data base is:\\
\textsf{http://www.research.att.com/~njas/sequences/index.html}};
$M(n)=\#\{ij|(i,j)\in\mathbb{N}^2_n\}$ in order to get some lower
bounds for it. To do this, we introduce the concept of
\textit{Universal Multiplication Table}, which is an infinite
array generated by multiplying the components of points in the
infinite lattice $\mathbb{N}^2$. Let $D(n)=\left\{d:d>0,~d|n
\right\}$. To get above mentioned bounds for the function $M(n)$,
we will need some upper bounds for the \textit{Divisor Function}
$d(n)=\#D(n)$, which we recall best known, due to J.L. Nicolas
\cite{jln}:
$$
\frac{\log d(n)}{\log 2}\leq\frac{\log n}{\log\log
n}\left(1+\frac{1.9349\cdots}{\log\log n}\right)\hspace{10mm}(n\geq
3),
$$
or
\begin{equation}\label{dn-nn}
d(n)\leq\mathfrak{N}(n)
\end{equation}
for $n\geq 3$, with
$$
\mathfrak{N}(n)=n^{\frac{\log 2}{\log\log
n}\left(1+\frac{387}{200\log\log n}\right)}.
$$

\section{Some Elementary Properties of the Function $\mathfrak{M}(n;k)$}
Considering (\ref{mnk}), for every $s\in\mathbb{C}$, we have:
\begin{equation}\label{1/ij-sum-mnk}
\sum_{1\leq i,j\leq
n}\frac{1}{(ij)^s}=\sum_{k=1}^{n^2}\frac{\mathfrak{M}(n;k)}{k^s}=
\sum_{k=1}^{\infty}\frac{\mathfrak{M}(n;k)}{k^s}.
\end{equation}
The left hand side of above identity is equal to $\zeta^2_n(s)$, in
which $\zeta_n(s)=\sum_{i=1}^n\frac{1}{i^s}$, and the number of
summands in the right hand side of above identity, is equal to
$M(n)$. Also, summing and counting all numbers in $MT_{n\times n}$,
we obtain respectively:
$$
\sum_{k=1}^{n^2}k\mathfrak{M}(n;k)=\left(\frac{n(n+1)}{2}\right)^2,
$$
and
$$
\sum_{k=1}^{n^2}\mathfrak{M}(n;k)=n^2,
$$
which both of them are special cases of (\ref{1/ij-sum-mnk}) for
$s=-1$ and $s=0$, respectively. To have some formulas for the
function $\mathfrak{M}(n;k)$, we define \textit{Incomplete Divisor
Function} to be $d(k;x)=\# D(k)\cap [1,x]$. This function has some
properties, which we list some of them:\\
\textbf{1.} It is trivial that for every $x\geq 1$ we have:
$$
1\leq d(k;x)\leq\min\left\{x,d(k)\right\}.
$$
So, $d(k;x)=O(x)$ and naturally we ask: What is the exact order of
$d(k;x)$? The next property, maybe useful to find answer.\\
\textbf{2.} If we let
$D(k)=\left\{1=d_1,d_2,\cdots,d_{d(k)}=k\right\}$, then we have:
\begin{eqnarray*}
\int_1^k
d(k;x)dx&=&\sum_{i=1}^{d(k)-1}(d_{i+1}-d_i)i=\sum_{i=1}^{d(k)-1}(i+1)d_{i+1}-
id_i-\sum_{i=1}^{d(k)-1}d_{i+1}\\&=&d(k)d_{d(k)}-1d_1-\sum_{d|k,d>1}d=kd(k)-\sigma(k),
\end{eqnarray*}
where $\sigma(k)=\sum\limits_{a\in D(k)}a$, and we have the
following bound due to G. Robin \cite{robin}:
\begin{equation}\label{sigma-rn}
\sigma(n)<\mathfrak{R}(n)\hspace{10mm}(n\geq 3),
\end{equation}
with
$$
\mathfrak{R}(n)=e^{\gamma}n\log\log n+\frac{3241n}{5000\log\log n},
$$
where $\gamma\approx 0.5772156649$ is Euler's constant. Considering
(\ref{dn-nn}) and (\ref{sigma-rn}), we obtain the following
inequality for every $k\geq 3$:
$$
2k-\mathfrak{R}(k)<\int_1^k d(k;x)dx<k\mathfrak{N}(k)-k-1.
$$
In general, every knowledge about $d(k;x)$ is useful, because:
\begin{prop} For every positive integers $k$ and $n$, we have:
$$
\mathfrak{M}(n;k)=d(k;n)-d\Big(k;\frac{k}{n}\Big)+R(n;k),
$$
where
$$
R(n;k)=\Big\lfloor\frac{k}{n}\Big\rfloor-\Big\lfloor\frac{k-1}{n}\Big\rfloor=
\begin{cases}
1,& n\mid k, \\
0, & {\rm other~wise}.
\end{cases}
$$
\end{prop}
\begin{proof} Considering (\ref{mnk}), we have:
$$
\mathfrak{M}(n;k)=\#\left\{(a,b)\in\mathbb{N}_n^2~|~ab=k\right\}=\sum_{d|k,d\leq
n,\frac{k}{d}\leq n}1=\sum_{d|k,\frac{k}{n}\leq d\leq n}1.
$$
Applying the definition of $d(k;x)$, completes the proof.
\end{proof}

\section{Universal Multiplication Table Function}
We define the \textit{Universal Multiplication Table Function}
$\mathfrak{M}(k)$ to be the number of $k$'s, which appear in the
universal multiplication table.
\begin{prop}\label{mk-dk} For every positive integer $k$, we have:
$$
\mathfrak{M}(k)=d(k).
$$
\end{prop}
\begin{proof} Here we have two proofs:\\
\textit{Elementary Method.} Considering the definition of universal
multiplication table, we have:
$$
\mathfrak{M}(k)=\lim_{n\rightarrow\infty}\mathfrak{M}(n;k)=\lim_{n\rightarrow\infty}\sum_{d|k,\frac{k}{n}\leq
d\leq n}1=\sum_{d|k,0<d<\infty}1=d(k).
$$
\textit{Analytic Method.} Considering (\ref{1/ij-sum-mnk}) for
$\Re(s)>1$ and taking limit both sides of it, when $n$ tends to
infinity, we obtain:
\begin{equation}\label{sum-mk/ks=zeta2}
\sum_{k=1}^{\infty}\frac{\mathfrak{M}(k)}{k^s}=\zeta^2(s),
\end{equation}
in which $\zeta(s)$ is the Riemann zeta-function. According to the
Theorem 11.17 of \cite{apos}, we obtain:
$$
\mathfrak{M}(k)=\lim_{T\rightarrow\infty}\frac{1}{2T}\int_{-T}^T\zeta^2(\sigma+it)k^{\sigma+it}dt,\hspace{10mm}(\sigma>1).
$$
Since $\zeta^2(s)=\sum_{m=1}^{\infty}d(m)m^{-s}$, we have:
\begin{eqnarray*}
\mathfrak{M}(k)&=&\lim_{T\rightarrow\infty}\frac{1}{2T}\int_{-T}^T\zeta^2(\sigma+it)k^{\sigma+it}dt\\
&=&\sum_{m=1}^{\infty}d(m)m^{-\sigma}k^{\sigma}\lim_{T\rightarrow\infty}\frac{1}{2T}\int_{-T}^T
\left(\frac{k}{m}\right)^{it}dt\\
&=&\sum_{m=1, m\neq
k}^{\infty}d(m)m^{-\sigma}k^{\sigma}\lim_{T\rightarrow\infty}
\frac{1}{2T}\int_{-T}^T\left(\frac{k}{m}\right)^{it}dt+d(k)\\
&=&\sum_{m=1, m\neq
k}^{\infty}d(m)m^{-\sigma}k^{\sigma}\lim_{T\rightarrow\infty}
\frac{1}{T}\sin\left(T\log\left(\frac{k}{m}\right)\right)+d(k)=d(k).
\end{eqnarray*}
This completes the proof.
\end{proof}
Now, fix positive integer $k$ and consider $\mathfrak{M}(n;k)$, as
an arithmetic function of the variable $n$. Clearly,
$\mathfrak{M}(n;k)$ is increasing, and for $n>k$, we have
$\mathfrak{M}(n;k)=\mathfrak{M}(k)$. Thus considering Proposition
\ref{mk-dk}, we obtain:
\begin{equation}\label{mnk-nk}
\mathfrak{M}(n;k)\leq d(k),
\end{equation}
and if $k\geq 3$, considering (\ref{dn-nn}) yields that:
$$
\mathfrak{M}(n;k)\leq \mathfrak{N}(k).
$$

\section{Statistical Study of $\mathfrak{M}(n;k)$'s}

Consider $S=\left[\mathfrak{M}(n;k)~|~1\leq k\leq n^2\right]$ as a
list of statistical data and suppose $\overline{\mathfrak{M}}(n)$ is
the average of above list, then we have:
$$
\overline{\mathfrak{M}}(n)=\frac{\sum_{k=1}^{n^2}\mathfrak{M}(n;k)}{\#\{ij|(i,j)\in\mathbb{N}^2_n\}}=
\frac{n^2}{M(n)}.
$$
Thus, we have:
\begin{equation}\label{mn-mn}
M(n)=\frac{n^2}{\overline{\mathfrak{M}}(n)}.
\end{equation}
Considering (\ref{mnk-nk}), it is clear that:
$$
\overline{\mathfrak{M}}(n)\leq\max\{\mathfrak{M}(n;k)\}_{k=1}^{n^2}\leq\max\{d(k)\}_{k=1}^{n^2}.
$$
To use (\ref{dn-nn}), we observe that the function $\mathfrak{N}(n)$
is increasing for $n\geq 114$. So, we have:
$$
\overline{\mathfrak{M}}(n)\leq\max\{d(1),d(2),\cdots,d(114),d(n^2)\}
\leq\max\{12,\mathfrak{N}(n^2)\}\hspace{10mm}(n\geq \sqrt{3}),
$$
and since $\mathfrak{N}(n)>114.1$ holds for every $n>0$, we obtain:
$$
\overline{\mathfrak{M}}(n)\leq \mathfrak{N}(n^2)\hspace{10mm}(n\geq
2).
$$
Therefore, we have proved the following result.
\begin{thm}\label{mnlb} For every $n\geq 2$, we have:
$$
M(n)\geq\frac{n^2}{\mathfrak{N}(n^2)}.
$$
\end{thm}
\begin{rem}\label{ford-comments}
One of the wonderful results about $MT_{n\times n}$ is
\textit{Erd\"{o}s Multiplication Table Theorem} \cite{pomerance},
which asserts:
$$
\lim_{n\rightarrow\infty}\frac{M(n)}{n^2}=0.
$$
Above theorem yields that in the Erd\"{o}s's theorem, however the
ratio $\frac{M(n)}{n^2}$ tends to zero, but it doesn't faster than
$\frac{1}{\mathfrak{N}(n^2)}$. More precisely, Erd\"{o}s showed
that $M(n)=n^2 (\log n)^{-c+o(1)}$ for $c=1+\frac{\log\log 2}{\log
2}$ \cite{erdos,ford}. The following table includes some
computational results about $M(n)$ by the Maple software.
\begin{center}
\begin{tabular}{||c|c|c||c|c|c||}
   \hline
  $n$   & $M(n)$  & $M(n)/n^2\approx$ & $n$   & $M(n)$  & $M(n)/n^2\approx$ \\ \hline
  10    & 42  & 0.4200000000 & 2000  & 959759 & 0.2399397500\\ \hline
  50    & 800  & 0.3200000000 &  3000  & 2121063 & 0.2356736667 \\ \hline
  100   & 2906  & 0.2906000000 & 4000  & 3723723 & 0.2327326875 \\ \hline
  1000  & 248083 & 0.2480830000 & 5000  & 5770205 & 0.2308082000 \\ \hline
\end{tabular}
\end{center}
Note that, the true order of $M(n)$ is $n^2 (\log n)^{-c}
(\log\log n)^{-3/2}$ \cite{ford}.
\end{rem}

\subsection*{Acknowledgment}I would like to express my gratitude
to Professor Aleksandar Ivic for his kind help on calculating
integral of $\mathfrak{M}(k)$. Also, I deem my duty to thank
Professor Jean-Louis Nicolas for kind sending the paper
\cite{jln}. Finally, I would like to thank professor Kevin Ford
for his very kind helps to clarify the historical background of
this note.

\end{document}